\documentclass[reqno0,12pt]{amsart}
\pdfoutput=1
\usepackage{dsfont}
\usepackage{bbm}
\usepackage[nodate]{datetime}
\usepackage[hmargin=27mm,top=30mm,bottom=30mm,a4paper]{geometry}
\usepackage{color}
\usepackage{subfig}
\usepackage{fancyhdr}
\usepackage{amsmath,amssymb,amsfonts,amsthm}
\usepackage{amstext}
\usepackage{amsmath}
\usepackage{amssymb}
\usepackage{enumerate}
\usepackage{amsbsy}
\usepackage{amsopn}
\usepackage{bbm,amsthm}
\usepackage{amscd}
\usepackage{amsxtra}
\usepackage{todonotes}
\usepackage{hyperref}

\newtheorem{theorem}{Theorem}[section]
\newtheorem*{theorem*}{Theorem}

\newtheorem{corollary}{Corollary}[section]

\theoremstyle{remark}

\newcommand{\CC}{\mathds{C}}

\newcommand{\RR}{\mathds{R}}
\newcommand{\NN}{\mathds{N}}
\newcommand{\HH}{\mathds{H}}

\newcommand{\ind}{\mathds{1}}

\newcommand{\EE}{\mathbb{E}}
\newcommand{\PP}{\mathbb{P}}

\begin{document}
\title{Random multiplicative functions: The Selberg-Delange class }
\author{Marco Aymone}
\begin{abstract}
Let $1/2\leq\beta<1$, $p$ be a generic prime number and $f_\beta$ be a random multiplicative function supported on the squarefree integers such that $(f_\beta(p))_{p}$ is an i.i.d. sequence of random variables with distribution $\mathbb{P}(f(p)=-1)=\beta=1-\mathbb{P}(f(p)=+1)$. Let $F_\beta$ be the Dirichlet series of $f_\beta$. We prove a formula involving measure-preserving transformations that relates the Riemann $\zeta$ function with the Dirichlet series of $F_\beta$, for certain values of $\beta$, and give an application. Further, we prove that the Riemann hypothesis is connected with the mean behavior of a certain weighted partial sums of $f_\beta$.
\end{abstract}
\maketitle

\section{Introduction.}  
We say that $f:\NN\to\CC$ is a multiplicative function if $f(nm)=f(n)f(m)$ for all non-negative integers $n$ and $m$ with $\gcd(n,m)=1$, and that $f$ has support on the squarefree integers if for any prime $p$ and any integer power $k\geq 2$, $f(p^k)=0$. An important example of such function is the M\"obius $\mu$: The multiplicative function supported on the squarefree integers such that at each prime $p$, $\mu(p)=-1$.

Many important problems in Analytic Number Theory can be rephrased in terms of the mean behavior of the partial sums of multiplicative functions. For instance, the Riemann hypothesis -- The statement that all the non-trivial zeros of the Riemann $\zeta$ function have real part equal to $1/2$ -- is equivalent to the statement that the partial sums of the M\"obius function have square root cancellation, that is,  $\sum_{n\leq x}\mu(n)$ is $O(x^{1/2+\epsilon})$, for all $\epsilon>0$. In this direction, the best result up to date is of the type $\sum_{n\leq x}\mu(n)=O(x\exp(-c(\log x)^\alpha))$, for some positive constant $c>0$, and some $0<\alpha<1$. Any improvement of the type $\sum_{n\leq x}\mu(n)=O(x^{1-\epsilon})$ for some $\epsilon>0$ would be a huge breakthrough in Analytic Number Theory, since it would imply that the Riemann $\zeta$ function has no zeros with real part greater than $1-\epsilon$. 

This equivalence between the Riemann hypothesis with the behavior of the partial sums of the M\"obius function led Wintner \cite{wintner} to investigate the behavior of a random model $f$ for the M\"obius function. This random model $f$ is defined as follows: We have that $f(n)$ is a random multiplicative function supported on the squarefree integers such that at primes $p\in\mathcal{P}$ (here $\mathcal{P}$ stands for the set of primes), $(f(p))_{p\in\mathcal{P}}$ is an i.i.d. sequence of random variables whith distribution $\PP(f(p)=-1)=\PP(f(p)=+1)=1/2$. It is important to obeserve that the sequence $(f(n))_{n\in\NN}$ is highly dependent, for instance, since $30=2\times 3\times 5$, we have that $f(30)$ depends on the values $f(2)$, $f(3)$ and $f(5)$. Wintner proved the square root cancellation for the partial sums of $f$, that is,
$\sum_{n\leq x}f(n)=O(x^{1/2+\epsilon})$ for all $\epsilon>0$, almost surely, and hence the assertion that \textit{the Riemann hypothesis is almost always true}. This upper bound has been improved several times: \cite{erdosuns}, \cite{halasz}, \cite{basquinn} and \cite{tenenbaum2013}. The best upper bound up to date is due to Lau, Tenenbaum and Wu \cite{tenenbaum2013}, which states that $\sum_{n\leq x}f(n)=O(\sqrt{x}(\log\log x)^{2+\epsilon})$ for all $\epsilon>0$, almost surely, and the best $\Omega$ result is due to Harper \cite{harpergaussian} which states that for any $A>5/2$, $\sum_{n\leq x}f(n)$ is not $O(\sqrt{x}(\log \log x)^{-A})$ almost surely.

Here we consider a slight different model for the M\"obius function. We start with a parameter $1/2\leq\beta\leq1$ and consider a random multiplictive function $f_\beta$ supported on the squarefree integers where at primes, $(f_\beta(p))_{p\in\mathcal{P}}$ is an i.i.d. sequence of random variables with $\PP(f_\beta(p)=-1)=\beta=1-\PP(f_\beta(p)=+1)$. For $\beta=1/2$ we recover the Wintner's model; For $\beta=1$, $f_1$ is the M\"obius $\mu$; And for $\beta<1$, we have that $f_\beta(n)$ equal to $\mu(n)$ with high probability as $\beta$ is taken to be close to $1$. In this paper we are interested in the following questions. 

\noindent \textit{Question 1}. What can be said about the partial sums $\sum_{n\leq x}f_\beta(n)$ for $1/2<\beta<1$? Does it have square root cancellation as in the Wintner's model and as we expect for the M\"obius function under the Riemann hypothesis?

\noindent \textit{Question 2}. If the partial sums $\sum_{n\leq x}f_\beta(n)$ are $O(x^{1-\delta})$ for some $\delta>0$, almost surely, then can we say something about the partial sums of the M\"obius function?

Considering the first question, observe that $\EE f(p)=1-2\beta$, and thus, we might say that at primes, $f_\beta(p)$ is equal to $1-2\beta$ in average. In the case $1/2<\beta<1$ the partial sums $\sum_{n\leq x}f_\beta(n)$ are well understood by the Selberg-Delange method, see the book of Tenenbaum \cite{tenenbaumlivro} chapter II.5 or the recent treatment of Granville and Koukoulopoulos \cite{koukouLSD}. Indeed, by the main result of \cite{koukouLSD}, we have that for $1/2<\beta<1$, the following holds almost surely 
$$\sum_{n\leq x}f_{\beta}(n)=(c_{f_\beta}+o(1))\frac{x}{(\log x)^{2\beta}},$$ where $c_{f_{\beta}}$ is a random constant which is positive almost surely. In particular, this implies that $\sum_{n\leq x}f_\beta(n)$ is not $O(x^{1-\delta})$, for any $\delta>0$, almost surely. This answers negatively our question 1. 

Here we provide a different proof of a negative answer to our question 1 for certain values of $\beta$, that is, the statement that we do not have square root cancellation for $\sum_{n\leq x}f_\beta(n)$ for certain values of $\beta$, almost surely. Further, by considering the question 2, we show that the Riemann hypothesis is equivalent to the square root cancellation of a certain weighted partial sums of $f_\beta$.

Before we state our results, let us introduce some notation. Given a probability space $(\Omega,\mathcal{F},\PP)$, let $\omega$ be a generic element of $\Omega$, and $T:\Omega\to\Omega$ be a measure-preserving transformation, \textit{i.e.}, $\PP(T^{-1}(A))=\PP(A)$, for all $A\in\mathcal{F}$. We see the random multiplicative function $f_\beta$ defined over the probability space $(\Omega,\mathcal{F},\PP)$ as a function $f_\beta:\NN\times\Omega\to\{-1,0,1\}$, that is, $f_\beta(n)$ is a random variable such that $f_\beta(n,\omega)\in\{-1,0,1\}$. Moreover, the Dirichlet series of $f_\beta$, say $F_\beta(s):=\sum_{n=1}^\infty\frac{f_\beta(n)}{n^s}$, is a random analytic function defined over the half plane $\HH_1:=\{s\in\CC:Re(s)>1\}$, that is $F_\beta:\HH_1\times \Omega\to\CC$ is such that $F_\beta(s,\omega)=\sum_{n=1}^\infty\frac{f_\beta(n,\omega)}{n^s}$ is analytic in the half plane $\HH_1$, for all $\omega\in\Omega$.

\begin{theorem}\label{teorema formula} Let  $n\geq 1$ be an integer, $\beta=1-\frac{1}{2^{n+1}}$, and $(\Omega,\mathcal{F},\PP)$ be a certain probability space, where it is defined $f_\beta$ for all values of $\beta\in[1/2,1]$. Let $F_\beta(s)=\sum_{n=1}^\infty\frac{f_\beta(n)}{n^s}$. Then there exists a measure-preserving transformation $T:\Omega\to\Omega$ such that $T^{2^n}=\mbox{identity}$ and such that the following formula holds for all $Re(s)>1$ and all $\omega\in\Omega$:
\begin{equation}\label{equation formula zeta}
\frac{1}{\zeta(s)^{2^n-1}}=\frac{1}{F_{1/2}(s,\omega)}\prod_{k=1}^{2^n} F_\beta(s,T^k\omega).
\end{equation}
\end{theorem}
In particular, if $\beta=3/4$, we have
\begin{equation*}
\frac{1}{\zeta(s)}=\frac{F_{3/4}(s,\omega)F_{3/4}(s,T\omega)}{F_{1/2}(s,\omega)}.
\end{equation*}
\begin{corollary}\label{corolario 1} For an integer $n\geq 1$ and $\beta=1-\frac{1}{2^{n+1}}$, we have that for any $\delta>0$, $\sum_{n\leq x}f_{\beta}(n)$ is not $O(x^{1-\delta})$ almost surely.
\end{corollary} 
Here we outline our proof of the Corollary \ref{corolario 1}. It utilizes the fact that the event in which $\sum_{n\leq x}f_\beta(n)=O(x^{1-\delta})$ is contained in the event in which the Dirichlet series $F_{\beta}(s)$ has analytic continuation to $\{Re(s)>1-\delta\}$. Moreover, one can easily check that for $\beta>1/2$, $F_\beta(1)=0$ almost surely. In the Wintner's proof \cite{wintner} of the square root cancellation of $\sum_{n\leq x} f_{1/2}(n)$, it has been proved that $F_{1/2}(s)$ is almost surely a non-vanishing analytic function over the half plane $\{Re(s)>1/2\}$. Thus, if we assume that for some $\delta>0$, $\sum_{n\leq x}f_\beta(n)=O(x^{1-\delta})$, almost surely, then the event  in which $F_\beta(s)$ has analytic continuation to $\{Re(s)>1-\delta\}$ also has probability $1$. Now the left side of \eqref{equation formula zeta} has a zero of multiplicity $2^n-1$ at $s=1$, and since $T$ preserves measeure, the right side of the same equation has a zero of multiplicity at least $2^n$ at the same point, which is a contradiction, and hence the event in which $F_\beta(s)$ has analytic continuation to $\{Re(s)>1-\delta\}$ can not hold with probability $1$. Moreover, by the Euler product formula for $Re(s)>1$ (here $\mathcal{P}$ stands for the set of primes)
\begin{equation}\label{equacao euler product}
F_\beta(s)=\prod_{p\in\mathcal{P}}\left(1+\frac{f_\beta(p)}{p^s} \right),
\end{equation}
 we see that the event in which $F_{\beta}$ has analytic continuation to $\{Re(s)>1-\delta\}$ is a tail event, in the sense that 
it does not depend in any outcome of a finite number of the random variables $f_{\beta}(p_1),...,f_{\beta}(p_r)$, where $p_1$,...,$p_r$ are primes. The Kolmogorov zero or one law states that each tail event has probability either equal to $0$ or $1$. Thus, the event in which $F_{\beta}$ has analytic continuation to $\{Re(s)>1-\delta\}$ has probability $0$, and hence the event in which $\sum_{n\leq x}f_\beta(n)=O(x^{1-\delta})$ also has probability $0$.

Now we turn our attention to Question 2. As mentioned above, clearly the event in which $\sum_{n\leq x}f_\beta(n)=O(x^{1-\delta})$ for some $\delta>0$ has probability $0$, and hence, the Question 2 as it is stated does not makes sense. However, if we consider a weighted sum of $f_\beta$, then  we can obtain an equivalence between the Riemann hypothesis with the mean behavior of a certain weighted partial sums of $f_\beta$. Before we state our next result, let $d(n)$ be the quantity of distinct primes that divide $n$.
\begin{theorem}\label{teorema equivalencia} The Riemann hypothesis is equivalent to 
$$\sum_{n\leq x}(2\beta-1)^{-d(n)}f_\beta(n)=O(x^{1/2+\epsilon}), \mbox{ for all }\epsilon>0,\mbox{ almost surely},$$ 
for each $\frac{1}{2}+\frac{1}{2\sqrt{2}}<\beta<1$.
\end{theorem}

\section{Preliminaries}
\subsection{Notation} Here we let $p$ denote a generic prime number and $\mathcal{P}$ to be the set of primes. We use $f(x)\ll g(x)$ and $f(x)=O(g(x))$ whenever there exists a constant $c>0$ such that $|f(x)|\leq c |g(x)|$, for all $x$ in a certain set $X$ -- This set $X$  could be all the interval $x\in[1,\infty)$ or $x\in (a-\delta,a+\delta)$, $a\in \RR,\delta>0$. We say that $f(x)=o(g(x))$ if $\lim\frac{f(x)}{g(x)}=0$. The notation $d|n$ means that $d$ divides $n$. Here $\ast$ stands for the Dirichlet convolution $(f\ast g)(n):=\sum_{d|n}f(d)g(n/d)$. We denote $d(n)=\sum_{p|n}1$, that is, the quantity of distinct primes that divides $n$. For a set $A$, $\ind_A(x)$ stands for the indicator function of the set $A$, that is, $\ind_A(x)=1$ if $x\in A$ and $\ind_A(x)=0$ if $x\notin A$.

\section{Proof of the results}
\subsection{Construction of the probability space}
We let $\mathcal{P}$ be the set of primes, $\Omega=[0,1]^{\mathcal{P}}=\{\omega=(\omega_p)_{p\in\mathcal{P}}:\omega_p\in[0,1]\mbox{ for all } p\}$, $\mathcal{F}$ the Borel sigma algebra of $\Omega$ and $\PP$ be the product of Lebesgue measures in $\mathcal{F}$. We set $f_\beta(p)$ as
$$f_\beta(p,\omega)=-\ind_{[0,\beta]}(\omega_p)+\ind_{(\beta,1]}(\omega_p).$$
It follows that $(f_\beta(p))_{p\in\mathcal{P}}$ are i.i.d. with distribution $\PP(f_\beta(p)=-1)=\beta=1-\PP(f_\beta(p)=+1)$. Also, we say that $f_\beta$ are uniformly coupled for different values of $\beta$. 
\subsection{Construction of the measure-preserving transformation} Now if $\beta=1-\frac{1}{2^{n+1}}$ with $n\geq 1$ an integer, we partionate the interval $[1/2,1]$ into $2^n$ subintervals $I_k=(a_{k-1},a_k]$ of lenght $\frac{1}{2^{n+1}}$ and with endpoints $a_k=\frac{1}{2}+\frac{k}{2^{n+1}}$. It follows that $a_0=1/2$, $a_{2^n-1}=\beta$ and $a_{2^n}=1$.  

Let $T_p:[0,1]\to[0,1]$ be the following interval exchange transformation: For $\omega_p\in[0,1/2]$, $T_p(\omega_p)=\omega_p$; In each interval $I_k$ as above the restriction $T_p|_{I_k}$ is a translation; $T_p(I_1)=I_{2^n}$ and for $k\geq 2$, $T_p(I_k)=I_{k-1}$. It follows that the $k$th iterate $T_p^k(I_k)=I_{2^n}$ and $T^{2^n}$ is the identity. Also, for each prime $p$, $T_p$ and its iterates preserve the Lebesgue measure and hence, $T:\Omega\to\Omega$ defined by $T\omega:=(T_p(\omega_p))_{p\in\mathcal{P}}$ preserves $\PP$, and so its iterates. 
\subsection{Proof of Theorem \ref{teorema formula}}
\begin{proof} We let $F_\beta$ be the Dirichlet series of $f_\beta$ and $I_k=(a_{k-1},a_k]$ be as above. Notice that $a_0=1/2$ and $a_{2^{n}}=1$, and hence $F_{a_0}=F_{1/2}$ and $F_{a_{2^n}}=F_1=\frac{1}{\zeta}$. Observe that
\begin{align*}
F_{1/2}\zeta=\frac{F_{a_0}}{F_{a_{2^n}}}=\frac{F_{a_0}}{F_{a_1}}\cdot \frac{F_{a_1}}{F_{a_2}}\cdot...\cdot \frac{F_{a_{2^n-1}}}{F_{a_{2^n}}}.
\end{align*}
Now, by the Euler product formula \eqref{equacao euler product}, we have that for all $Re(s)>1$
\begin{align*}
\frac{F_{a_k}}{F_{a_{k+1}}}(s,\omega)=\prod_{p\in\mathcal{P}}\frac{1+\frac{f_{a_k}(p,\omega_p)}{p^s}}{1+\frac{f_{a_{k+1}}(p,\omega_p)}{p^s}}=
\prod_{p\in\mathcal{P}} \frac{p^s+\ind_{I_k}(\omega_p)}{p^s-\ind_{I_k}(\omega_p)}.
\end{align*}
Thus, as all $I_k$ have same lenght, we see that each $\frac{F_{a_k}}{F_{a_{k+1}}}$ is equal in probability distribution to the last $\frac{F_{a_{2^n-1}}}{F_{a_{2^n}}}$. Moreover, if $T$ is as above, since $\ind_{I_k}(\omega_p)=\ind_{I_{2^n}}\circ T_p^k(\omega_p)$, we have that
$$\frac{F_{a_k}}{F_{a_{k+1}}}(s,\omega)=\frac{F_{a_{2^n-1}}}{F_{a_{2^n}}}(s,T^k\omega)=F_\beta(s,T^k\omega)\zeta(s).$$
Thus
$$F_{1/2}(s,\omega)\zeta(s)=\zeta(s)^{2^n}\prod_{k=1}^{2^n}F_\beta(s,T^k\omega),$$
which concludes the proof. \end{proof}
\subsection{Proof of Corollary \ref{corolario 1}}
\begin{proof}
A standard result about Dirichlet series, is that the Dirichlet series of an arithmetic function $f$, say $F(s)$, is the Mellin transform of the partial sums of $f$. Indeed, we have that for $s$ in the half plane of convergence of $F(s)$,
$$F(s)=s\int_1^\infty \frac{\sum_{n\leq x}f(n)}{x^{s+1}}dx.$$ 
Thus, we can conclude that the event in which the partial sums $\sum_{n\leq x}f(n)$ are $O(x^{\alpha})$ is contained in the event in which the Dirichlet series $F(s):=\sum_{n=1}^\infty \frac{f(n)}{n^s}$ is analytic in the half plane $\{Re(s)>\alpha\}$. Thus, under the assumption that $\sum_{n\leq x}f_\beta(n)=O(x^{1-\delta})$ almost surely, we have that $F_\beta(s)=\sum_{n=1}^\infty \frac{f_\beta(n)}{n^s}$ has analytic continuation to the half plane $\{Re(s)>1-\delta\}$ almost surely. Moreover, we can check that $F_\beta(1)=0$ almost surely. Indeed, by taking the logaritihm of the Euler product formula \eqref{equacao euler product} and then using Taylor expansion of each logarithm, we see that
\begin{equation}\label{equacao exponencial euler product}
F_\beta(s)=\exp\left(\sum_{p\in\mathcal{P}}\frac{f_\beta(p)}{p^s}+A_{\beta}(s) \right),
\end{equation}
where $A_{\beta}(s)=O_{\sigma_0}(1)$ for all $Re(s)\geq \sigma_0>1/2$. Since $\EE f_{\beta}(p)=1-2\beta<0$ for all primes $p$, we have by the Kolmogorv two series Theorem that $\lim_{s\to 1^+}\sum_{p\in\mathcal{P}}\frac{f_\beta(p)}{p^s}=-\infty$ almost surely, and hence, $\lim_{s\to 1^+}F_\beta(s)=0$ almost surely.

If $T$ is the meausre-preserving transformation as in Theorem \ref{teorema formula}, then the same is almost surely true for $F_\beta(s,T^k\omega)$. Further, in the Wintner's proof \cite{wintner} of the square root cancellation of $\sum_{n\leq x} f_{1/2}(n)$, it has been proved that $F_{1/2}(s)$ is almost surely a non-vanishing analytic function over the half plane $\{Re(s)>1/2\}$. Indeed, this can be proved by the formula \eqref{equacao exponencial euler product}. 

A well known fact is that the Riemann $\zeta$ function has a simple pole at $s=1$, and hence, $\frac{1}{\zeta(s)}$ has a simple zero at the same point. Moreover we recall that if an analytic function $G$ has a zero at $s=s_0$, then there exists a non-vanishing analytic function $H$ at $s=s_0$ and a non-negative integer $m$, called the multiplicity of the zero $s_0$, such that $G(s)=(s-s_0)^mH(s)$. Thus the left side of 
\begin{equation*}\label{equacao formula zeta}
\frac{1}{\zeta(s)^{2^n-1}}=\frac{1}{F_{1/2}(s,\omega)}\prod_{k=1}^{2^n} F_\beta(s,T^k\omega).
\end{equation*}
 has a zero of multiplicity $2^n-1$ at $s=1$ while the right side of the same equation has a zero of multiplicity at least $2^n$ at the same point, almost surely, which is a contradiction. Thus we see that the probability of the event in which $F_{\beta}(s)$ has analytic continuation to $Re(s)>1-\delta$ is strictly less than one. Now we can check by the Euler product formula \eqref{equacao euler product} that the event in which $F_\beta$ has analytic continuation to $Re(s)>1-\delta$ is a tail event for $\delta<1$, \textit{i.e.}, whether $F_\beta$ has analytic continuation to $\{Re(s)>1-\delta\}$ does not depend in any outcome of a finite number of random variables $\{f_{\beta}(p):p\leq y\}$. Indeed, we can write 
 \begin{equation*}
F_\beta(s)=\prod_{p\leq y}\left(1+\frac{f_\beta(p)}{p^s} \right)\prod_{p> y}\left(1+\frac{f_\beta(p)}{p^s} \right),
\end{equation*} 
and since $\prod_{p\leq y}\left(1+\frac{f_\beta(p)}{p^s} \right)$ is a non-vanishing analytic function in $Re(s)>0$, we obtain that $F_{\beta}(s)$ has analytic continuation to $Re(s)>1-\delta$ ($\delta<1$) if and only if $X_y(s):=\prod_{p> y}\left(1+\frac{f_\beta(p)}{p^s} \right)$ has analytic continuation to the same half plane. Since $X_y(s)$ is independent of  $\{f_{\beta}(p):p\leq y, p\in\mathcal{P}\}$ and the random variables $(f_\beta(p))_{p\in\mathcal{P}}$ are independent, we conclude that the event in which $F_\beta$ has analytic continuation to $\{Re(s)>1-\delta\}$ is a tail event.

Thus by the Kolmogorov zero or one law, we have that the probability in which $F_\beta$ has analytic continuation to $\{Re(s)>1-\delta\}$ is zero, and hence the probability of $\sum_{n\leq x}f_\beta(n)=O(x^{1-\delta})$ is also zero. 
\end{proof}

\subsection{Proof of Theorem \ref{teorema equivalencia}}
\begin{proof} We begin by observing that the function $g_\beta(n):=(2\beta-1)^{-d(n)}f_\beta(n)$ is multiplicative and supported on the squarefree integers. Moreover, at each prime $p$, $g_\beta(p)=\frac{f_\beta(p)}{2\beta-1}$, and hence $\EE g_\beta(p)=-1$. If $\beta>\frac{1}{2}+\frac{1}{2\sqrt{2}}$, we have that
$$A_\beta(s):=\sum_{p\in\mathcal{P}}\sum_{m=2}^{\infty}\frac{(-1)^{m+1}}{m}\frac{g_\beta(p)^m}{p^{ms}}$$
converges absolutely for all $Re(s)>1/2$, and hence defines a random analytic function in this half plane. Moreover, $A_\beta(s)=O_{\sigma_0}(1)$ uniformly for all $Re(s)\geq \sigma_0>1/2$. Thus, by the Euler product formula \eqref{equacao euler product} for $g_\beta$, we have that the Dirichlet series $G_{\beta}(s):=\sum_{n=1}^\infty\frac{g_\beta(n)}{n^s}$ can be represented in the half plane $Re(s)>1$ as
$$G_{\beta}(s)=\exp\left(\sum_{p\in\mathcal{P}}\frac{g_\beta(p)}{p^s}+A_\beta(s) \right).$$
Moreover, by the same argumet, there exists an analytic function $B(s)$ with the same properties of $A_\beta(s)$ such that
$$\zeta(s)=\exp\left(\sum_{p\in\mathcal{P}}\frac{1}{p^s}+A_\beta(s) \right).$$
Now observe that
$$H_\beta(s):=G_{\beta}(s)\zeta(s)=\exp\left(\sum_{p\in\mathcal{P}}\frac{g_\beta(p)+1}{p^s}+A_\beta(s)+B(s) \right).$$
Now, by the Kolmogorov one series Theorem, $\sum_{p\in\mathcal{P}}\frac{g_\beta(p)+1}{p^s}$ converges almost surely for all $Re(s)>1/2$, and hence it defines, almost surely, a random analytic function in this half plane. Moreover, by Theorem 3.1 of \cite{aymonebiased}, for fixed $1/2<\sigma\leq 1$, we have that for all large $t>0$, $\sum_{p\in\mathcal{P}}\frac{g_\beta(p)+1}{p^{\sigma+it}}\ll (\log t)^{1-\sigma}\log\log t$, almost surely. Thus, for each fixed $1/2<\sigma$, 
$$H_\beta(\sigma+it),1/H_\beta(\sigma+it)\ll t^{\epsilon},$$ for all $\epsilon>0$, almost surely. A well known consequence of the Riemann Hypothesis, is that $1/\zeta(s)$ has analytic continuation to $Re(s)>1/2$ and for each fixed $\sigma>1/2$, $1/\zeta(\sigma+it)\ll t^\epsilon$, for all $\epsilon>0$. Thus, if we assume the Riemann hypothesis, we obtain that $G_{\beta}(s)$ has analytic continuation to $Re(s)>1/2$ given by $G_\beta(s)=H_\beta(s)/\zeta(s)$, and for each fixed $\sigma>1/2$, $G_\beta(\sigma+it)\ll t^\epsilon$ for all $\epsilon>0$, almost surely. Now a convergence result for Dirichlet series (see for instance Theorem 2.8, page 223 of \cite{tenenbaumlivro}) gives that $G_\beta(s)$ converges for all $Re(s)>1/2$. Now by Kroenecker's Lemma (see \cite{shiryaev} page 390), we have that $\sum_{n\leq x}g_\beta(n)\ll x^{1/2+\epsilon}$ for all $\epsilon>0$, almost surely. On the other hand, if $\sum_{n\leq x}g_\beta(n)\ll x^{1/2+\epsilon}$ for all $\epsilon>0$, almost surely, then $G_\beta(s)$ is almost surely analytic in $Re(s)>1/2$, and thus $G_\beta(s)/H_\beta(s)$ also is almost surely analytic in $Re(s)>1/2$. Since, $1/\zeta(s)=G_\beta(s)/H_\beta(s)$, we have that $1/\zeta(s)$ has analytic continuation to $Re(s)>1/2$. This last assertion is equivalent to the Riemann hypothesis.   
\end{proof}

Departamento de Matem\'atica, Universidade Federal de Minas Gerais, Av. Ant\^onio Carlos, 6627, CEP 31270-901, Belo Horizonte, MG, Brazil. \\
\textit{Email address:} aymone.marco@gmail.com

\end{document}